\documentclass[12pt,reqno]{article}
\usepackage{amssymb}
\usepackage{amsmath}
\usepackage{mathrsfs}
\usepackage{amsthm}
\usepackage{amscd}
\usepackage{amsbsy}
\usepackage{amsfonts}
\usepackage{appendix}
\usepackage{float}
\usepackage{xcolor}
\usepackage{epstopdf}
\usepackage{graphicx}
\usepackage{psfrag}
\usepackage{subfigure}
\usepackage{overpic}
\usepackage{longtable}
\usepackage{multirow}
\usepackage{verbatim}
\usepackage{listings}
\usepackage{epsfig}
\usepackage{enumerate}
\usepackage{booktabs}
\usepackage{bm}
\usepackage{stackengine}
\newtheorem{thm}{Theorem}[section]

\newtheorem{lem}[thm]{Lemma}
\newtheorem{prop}[thm]{Proposition}

\newtheorem{rem}[thm]{Remark}

\newcommand\xrowht[2][0]{\addstackgap[.5\dimexpr#2\relax]{\vphantom{#1}}}
\renewcommand{\theequation}{\thesection.\arabic{equation}}

\def\R{{\mathbb R}}

\def\bc{\begin{center}}       \def\ec{\end{center}}
\def\be{\begin{equation}}     \def\ee{\end{equation}}
\def\ba{\begin{array}}        \def\ea{\end{array}}
\def\bea{\begin{eqnarray}}    \def\eea{\end{eqnarray}}
\def\beaa{\begin{eqnarray*}}  \def\eeaa{\end{eqnarray*}}

\def\la{\lambda}

\def\oo{\infty}               

\def\dd{\cdots}               
                 \def\qq{\qquad}

\def\Lp{{\mathcal L}^p}

\def\lb{\label}
\def\x#1{{\rm (\ref{#1})}}

\def\ifl{\iffalse}
\def\Proof{\noindent{\bf Proof} \quad}
\def\qed{\hfill $\Box$ \smallskip}

\begin{document}

\title{\bf Non-integrability of the Critical Systems for Optimal Sums of Eigenvalues of Sturm-Liouville Operators}

\author{Yuzhou Tian$^{a,}$\footnote{Correspondence author.},\qq Meirong Zhang$^{b}$}

\date{\today}%

\maketitle

\begin{center}
$^a$ Department of Mathematics, Jinan University, Guangzhou 510632, China	\\
$^b$ Department of Mathematical Sciences, Tsinghua University, Beijing 100084, China\\
E-mail: {\tt
tianyuzhou2016@163.com (Y. Tian)}\\
E-mail: {\tt zhangmr@tsinghua.edu.cn (M. Zhang)}
\end{center}

\begin{abstract}
 The optimal lower or upper bounds for sums of the first $m$ eigenvalues of Sturm-Liouville operators can be obtained by solving the corresponding critical systems, which are Hamiltonian systems of $m$ degrees of freedom with $m$ parameters. With the help of the differential Galois theory, we prove that  these critical systems are not meromorphic integrable except for two known completely  integrable cases. The non-integrability of the critical systems reveal  certain  complexities for the original eigenvalues problems.
\end{abstract}

{\small {\bf Mathematics Subject Classification (2020)}:
Primary 34L15; 
Secondary 37J30, 
70H07. 

{\bf Keywords:} Sums of eigenvalues, Sturm-Liouville operator, Critical system, Meromorphic integrable, Differential Galois theory.}


\section{Introduction} \lb{mr}
We consider the following Dirichlet eigenvalue problem for the Sturm-Liouville operator
\begin{equation}\label{line}
	-y'' + q \left(x\right)y=\la y, \qq x\in\left[0,1\right],
\end{equation}
subject to the Dirichlet boundary condition
\begin{align}\label{dc}
&y\left(0\right)=y\left(1\right)=0,
\end{align}
where $q\left(x\right)\in\Lp:=L^p\left(\left[0,1\right],\R\right)$ is 
an integrable potential with the exponent $p\in\left(1,\infty\right)$. The real number $\la$ is called an eigenvalue of system \eqref{line} if there is a nontrivial solution $y\left(x\right)$ of system \eqref{line}. It is
well-known that problem  \eqref{line}-\eqref{dc} exists
a sequence of eigenvalues
 \[
 \la_1(q) < \la_2(q) < \dd < \la_m(q) < \dd, \qq \la_m(q) \to +\oo \mbox{ as } m\to \oo.
 \]
In the spectral theory of Sturm-Liouville operator, the optimal estimates of eigenvalues have been the focus of investigation for mathematicians and physicists due to the that such eigenvalues problems are very useful for understanding the nonhomogeneous string equation  and various dimensionally reduced  wave equations, see \cite{MR4339006,MR73776}, etc.  Since Krein's seminal work in \cite{MR73776},  researchers have done substantive work and achieved enormous progress for the optimal estimates of eigenvalues.  In the last few decades, the important developments of eigenvalues problems for Sturm-Liouville operator \eqref{line} mainly involve gaps \cite{MR1081670,MR2881964,MR4447102,MR1948113,TWZ}, individual eigenvalues \cite{MR4630844,MR4673607}, ratios \cite{MR1218744, MR4700368, MR4339006}, nodes \cite{chu2025,chu2024,MR4444221} and  applications in nonlinear PDEs \cite{MR899179}.

Besides the above  eigenvalues problems, the estimates for sum of eigenvalues of ordinary differential operators are also very important in the spectral theory, see \cite{MR1744872,MR4252029,MR3392903,MR75382,MR77748} and the references therein. For Sturm-Liouville operator \eqref{line}, the sum of the first $m$ eigenvalues can be described as the following optimization problems
\begin{align}\label{eq66}
	&\mathscr{E}_m^-:=\min\limits_{q\in B_{p,r}}\sum_{i=1}^m\lambda_i\left(q\right)\;\text{and}\; \mathscr{E}_m^+:=\max\limits_{q\in B_{p,r}}\sum_{i=1}^m\lambda_i\left(q\right),
\end{align}
where 
$B_{p,r}=\left\{q\in \Lp:\;\parallel q\parallel_p\leq r\right\}$
is a ball of the Lebesgue space $\left(\Lp,\parallel\cdot\parallel_p\right)$. Very recently,
  an equivalent characterization to the optimization problems \eqref{eq66} was shown by Tian and Zhang \cite{TZ}, who presented a general method to obtain the optimal lower or upper bound for sum of the first $m$  eigenvalues of \eqref{line} with the exponent $p\in(1,+\oo)$. More precisely, the optimal potentials of  problems \eqref{eq66} can be determined by the solvability of the critical systems consisting of nonlinear ordinary differential equations, as shown in Theorem 2.5 of \cite{TZ}.
  
  With suitable selection of exponent $p=k/\left(k-1\right)$, $k\in\mathbb{N}^+$, the critical system is equivalent to the following polynomial Hamiltonian systems of $m$ degrees of freedom
  \begin{align}\label{eq38}
  	&u'_i=v_i,\quad v'_i=-\mu_iu_i+\varepsilon\left(\sum_{j=1}^mu_j^2\right)^{k-1}u_i,\quad  i=1,\ldots,m,
  \end{align}
  with the Hamiltonian
  \begin{align}\label{eq39}
  	&H=\frac{1}{2}\sum_{i=1}^m\left(v_i^2+\mu_iu_i^2\right)-\frac{\varepsilon}{2k}\left(\sum_{j=1}^mu_j^2\right)^k
  \end{align}
and $\varepsilon=\pm$, which was given by (3.19) and (3.20) of \cite{TZ}. The paper \cite{TZ} has shown that the integrability of Hamiltonian system \eqref{eq38} is of great benefit to  the solvability of the critical system. However, there are no universal techniques  to decide the integrability of Hamiltonian systems. To the best of our knowledge, most studies have been dedicated to the integrability for Hamiltonian system of  two degrees of freedom, see for instance \cite{MR3743739,MR4075569,MR2831794,MR943701} and references therein. Except the Hamiltonian system with homogeneous potential, there are a few papers relating to the integrability of other types of Hamiltonian systems of arbitrary degrees of freedom, see \cite{MR1764944,MR2123446,MR3804721}. Tian and Zhang  found two completely integrable cases of Hamiltonian system \eqref{eq38}: $k=2$ and all parameters $\left(\mu_1,\ldots,\mu_m\right)\in\mathbb{R}^m$; all integers $k\geq3$ and $\mu_1=\mu_2=\cdots=\mu_m$, see Theorem 3.4 of \cite{TZ}. In addition, based on the numerical simulations, they conjectured that system \x{eq38} is non-integrable except for the above two cases.

Our present paper will focus on  the integrability of Hamiltonian system \eqref{eq38} in sense of meromorphic. We will show that system \eqref{eq38}  is not meromorphic integrable except for the integrable cases as pointed in  \cite{TZ} and mentioned above. See Theorem \ref{th1} of Section \ref{se5}. Our result gives a positive answer to the  conjecture of \cite{TWZ}. We employ the differential Galois theory to analyse the meromorphic non-integrability of Hamiltonian system \eqref{eq38}. Although this technique on the theoretical level is mature, it is not easy to use it to practical problems because there is no general approach to analyse the structures of the corresponding differential Galois groups. Compared with the previous results \cite{MR3743739,MR4075569,MR2831794,MR943701}, we need to analyse the variational equations along two particular solutions rather than one, which will clearly demonstrate structures of their differential Galois groups.

The structure of the paper is as follows. In Section \ref{se5}, we first recall a necessary condition for
meromorphic integrability of Hamiltonian systems, which is known as Morales-Ramis theory. Then, we devote to prove the meromorphic non-integrability of Hamiltonian system \eqref{eq38}, see Lemmas \ref{le2} and \ref{le3}. Finally, we give a complete classification of
meromorphic integrability for Hamiltonian system \eqref{eq38}, see Theorem \ref{th1}. The classic hypergeometric differential equation  \cite{MR277789} and Kovacic's algorithm \cite{MR839134} will play an important roles in our proof. For convenience, the hypergeometric differential equation and Kovacic's algorithm were collected in Appendixes \ref{A} and \ref{B}.
\section{Non-integrability of system \eqref{eq38} }\label{se5}
The main topics of this section will be investigated  the non-integrability of Hamiltonian system \eqref{eq38} in sense of meromorphic. 

Consider a complex symplectic manifold $M\subset\mathbb{C}^{2m}$ of dimension $2m$  with the standard symplectic form $\bm{\tilde{\omega}}=\sum_{j=1}^mdu_j\wedge dv_j$.  Let $H: M\rightarrow \mathbb{C}$ be a holomorphic  Hamiltonian. The Hamiltonian system with $m$ degrees of freedom is given by
\begin{align}\label{eq59}
	&\dfrac{d \mathbf{x}}{dt}=X_H\left(\mathbf{x}\right)=\left(\dfrac{\partial H}{\partial \mathbf{v}}, -\dfrac{\partial H}{\partial \mathbf{u}}\right),\quad t\in\mathbb{C}, \quad \mathbf{x}=\left(\mathbf{u},\mathbf{v}\right)\in M,
\end{align}
where $\mathbf{u}=\left(u_1,\ldots,u_m\right)$ and $\mathbf{v}=\left(v_1,\ldots,v_m\right)$ are the canonical coordinates. We say that Hamiltonian system \eqref{eq59} is completely integrable or Liouville
integrable if there exists $m$ non-constant functions $I_1\equiv H,I_2,\ldots,I_m$ satisfying the following conditions.
\begin{itemize}
	\item [(i)] The functions $I_i$ for $i=1,\ldots,m$ are functionally independent, that is, their gradients $\nabla I_1,\ldots,\nabla I_m$ are
linearly independent.
		\item [(ii)] The Poisson bracket of $I_i$ and $I_j$ is 
$$\left\{I_i,I_j\right\}=\sum_{l=1}^m\left(\frac{\partial I_i}{\partial v_l}\frac{\partial I_j}{\partial u_l}-\frac{\partial I_i}{\partial u_l}\frac{\partial I_j}{\partial v_l}\right)=0$$
for all $i,j=1,\dots,m$.
\end{itemize}
The above functions $I_i$ for $i=1,\ldots,m$ are called first integrals of system \eqref{eq59}. In addition, Hamiltonian system \eqref{eq59} is meromorphic completely integrable, or simply meromorphic integrable if its $m$ functionally independent first integrals $I_1\equiv H,I_2,\ldots,I_m$ are meromorphic. 

Let $\Gamma$  be  a non-equilibrium solution to system  \eqref{eq59}. Assume that $\Gamma$ can be parameterized by time $t$, that is,
\begin{align*}
	\bm{\varphi}:\mathbb{C}&\rightarrow M\subset\mathbb{C}^{2m}\\
	t&\mapsto\left(\mathbf{u}\left(t\right),\mathbf{v}\left(t\right)\right).
\end{align*}
Then the variational equation (VE,  for short) along $\Gamma$ is the linear differential system
\begin{align}\label{eq60}
	&\dfrac{d \mathbf{y}}{dt}=\dfrac{\partial X_H\left(\bm{\varphi}\left(t\right)\right)}{\partial \mathbf{x}}\mathbf{y}, \quad \mathbf{y}\in T_\Gamma M,
\end{align}
where $T_\Gamma M$ is the tangent bundle $TM$ restricted on $\Gamma$.

Let $N:=T_\Gamma M/T\Gamma$ be the normal bundle
of $\Gamma$ \cite{MR1411677}, and $\pi: T_\Gamma M \rightarrow N$ be the nature projective homomorphism. The normal variational equation (NVE,  for short) along $\Gamma$ has the form
\begin{align}\label{eq61}
	&\dfrac{d \mathbf{z}}{dt}=\pi_*\left(T\left(\mathfrak{u}\right)\left(\pi^{-1}\mathbf{z}\right)\right),\quad \mathbf{z}\in N,
\end{align}
where $\mathfrak{u}=X_H\left(\mathbf{x}\right)$ with $\mathbf{x}\in M$, and $T\left(\mathfrak{u}\right) $ is the tangential
variation of $\mathfrak{u}$ along $\Gamma$, that is, $T\left(\mathfrak{u}\right)=\partial X_H/\partial \mathbf{x}$.  Note that the above NVE is a $2\left(m-1\right)$-dimensional linear differential system. We can employ a generalization of D'Alambert's method to get the NVE \eqref{eq61}, see \cite{MR1713573}.
Briefly speaking, we use the fact that $X_H\left(\bm{\varphi}\left(t\right)\right)$ is a solution of the VE \eqref{eq60} to reduce its dimension  by one. In effect, we typically restrict the equation \eqref{eq59} to the energy level $h=H\left(\bm{\varphi}\left(t\right)\right)$. Then the dimensionality of  the corresponding VE \eqref{eq60} also can be reduced.

Morales and Ramis \cite{MR1713573} proved the following classical theorem, which gives a necessary condition for the integrability of Hamiltonian system \eqref{eq59} in the Liouville sense.
For precise notions of differential Galois theory, see \cite{MR1960772}.
\begin{thm}\label{th5}{\rm (Morales-Ramis theorem, see \cite{MR1713573})} If Hamiltonian system \eqref{eq59} is meromorphically integrable in the Liouville
	sense in a neighbourhood of a particular solution $\Gamma$, then the identity component of the Galois group of the NVE \eqref{eq61} is Abelian.
\end{thm}
Theorem \ref{th5} establishes a relation between the meromorphic  integrability and the differential Galois group of the NVE \eqref{eq61}. However, it's very difficult to decide the solvability of the differential Galois group of NVE \eqref{eq61} in applications.

The next theorem tells us that the identity component of the differential Galois group is invariant under the covering.
\begin{thm}\label{th4}{\rm (\cite{MR1713573})}
	Let $\mathcal{M}$ be a connected Riemann surface and $\nabla$ be a meromorphic connection over $\mathcal{M}$. 	Assume that $f:\mathcal{M}'\longrightarrow\mathcal{M}$ is a finite ramified covering of $\mathcal{M}$ by a connected Riemann surface  $\mathcal{M}'$.  Let $\nabla'=f^*\nabla$ be the pull back of $\nabla$ by $f$. Then there exists a natural
	injective homomorphism
	$$\emph{Gal}\left(\nabla'\right)\rightarrow\emph{Gal}\left(\nabla\right)$$
	of differential Galois groups which induces an isomorphism between their Lie
	algebras.
\end{thm}

Next, we discuss the meromorphic non-integrability of Hamiltonian system \eqref{eq38}.
\begin{lem}\label{le2}
If $k\geq3$, $\mu_1\neq\mu_2$ and $\mu_1\mu_2\neq0$, then Hamiltonian system \eqref{eq38}  is not meromorphic integrable.
\end{lem}

\Proof
One can easily observe that system \eqref{eq38} has two invariant manifolds
\begin{align*}
&\mathcal{N}_1=\left\{\left(\mathbf{u},\mathbf{v}\right)\in \mathbb{C}^{2m}\mid u_j=v_j=0,j=2,\ldots,m\right\},\\
&\mathcal{N}_2=\left\{\left(\mathbf{u},\mathbf{v}\right)\in \mathbb{C}^{2m} \mid u_1=v_1=0, u_j=v_j=0,j=3,\ldots,m\right\}.
\end{align*}
System \eqref{eq38} restricted to the first invariant manifold $\mathcal{N}_1$  becomes
\begin{align}\label{eq1}
	&u'_1=v_1,\quad v'_1=-\mu_1u_1+\varepsilon u_1^{2k-1},
\end{align}
with the Hamiltonian
\begin{align}\label{eq16}
&h=\dfrac{1}{2}v_1^2+\dfrac{1}{2}\mu_1u_1^2-\dfrac{\varepsilon}{2k}u_1^{2k}.
\end{align}
Solving equation \eqref{eq16}, one can get
\begin{align}\label{eq47}
&\dfrac{du_1}{dt}=\pm\sqrt{2h+\dfrac{\varepsilon}{k}u_1^{2k}-\mu_1u_1^2}.
\end{align}
Integrating equation \eqref{eq47}, we have 
\begin{align}\label{ceq47}
&t= \pm\displaystyle\int\frac{du_1}{\sqrt{2h+\dfrac{\varepsilon}{k}u_1^{2k}-\mu_1u_1^2}}.
\end{align}
For $k\geq3$, equation \eqref{ceq47} is called  hyperelliptic integrals, which are not always in terms of elementary functions, see \cite{MR0060642}.

Let $\Theta\left(h\right)\in\mathbb{C}^2$ be a solution  of system \eqref{eq1} on the energy level $h$. Consequently,
\begin{align}\label{eq57}
&\varGamma_h:=\left\{\left(u_1\left(t\right),v_1\left(t\right),0,\ldots,0\right)\in\mathbb{C}^{2m}\mid \left(u_1\left(t\right),v_1\left(t\right)\right)\in\Theta\left(h\right)\right\}
\end{align}
is a particular solution of system \eqref{eq38}. To perform Theorem \ref{th5}, we must find a non-equilibrium  particular solution $\varGamma_h$. Equation \eqref{eq47} has three equilibrium points  $u_1=0$ and $u_1=\pm\sqrt[2k-2]{\mu_1k/\varepsilon}$ on the energy level  $h=0$.  After eliminating these equilibria, we can construct a non-constant
solution $u_1\left(t\right)$, which provides a non-equilibrium  particular solution $\Gamma_0\in\varGamma_0$.

Let
$\bm{\xi}:=\left(\xi_1,\ldots,\xi_m\right)^T$ and  $\bm{\tilde{\xi}}:=\left(\tilde{\xi}_1,\ldots,\tilde{\xi}_m\right)^T$.  We obtain that the variational equation (VE) along $\Gamma_0$ is
\begin{align}\label{eq2}
	&\left(
	\begin{array}{c}
		\bm{\xi}'\vspace{1ex}\\
		\bm{\tilde{\xi}}'\\
	\end{array}\right)
	=
	\left(
	\begin{array}{cc}
		\mathbf{0}&\mathbf{I}\vspace{2ex}\\
		\mathbf{\Lambda}&\mathbf{0}\\
	\end{array}
	\right)
	\left(
	\begin{array}{c}
		\bm{\xi}\vspace{1ex}\\
		\bm{\tilde{\xi}}\\
	\end{array}\right),
\end{align}
where
$$\mathbf{\Lambda}:=\text{diag}\;\Big(\varepsilon\left(2k-1\right)u_1^{2k-2}\left(t\right)-\mu_1,\varepsilon u_1^{2k-2}\left(t\right)-\mu_2,\varepsilon u_1^{2k-2}\left(t\right)-\mu_3,\ldots,\varepsilon u_1^{2k-2}\left(t\right)-\mu_m\Big).$$
The VE \eqref{eq2} is composed of  $m$ uncoupled Schr\"{o}dinger equations
\begin{align*}
	&\bm{\xi}''=\mathbf{\Lambda}\bm{\xi},
\end{align*}
that is,
\begin{align}
	&\xi''_1=\left(\varepsilon \left(2k-1\right)u_1^{2k-2}\left(t\right)-\mu_1\right)\xi_1,\label{eq62}\\
	&\xi''_j=\left(\varepsilon u_1^{2k-2}\left(t\right)-\mu_j\right)\xi_j,\quad j=2,\ldots,m.
\end{align}
Since $\xi_1=u'_1\left(t\right)$ is a solution of \eqref{eq62},  equation  \eqref{eq62} can be solved by Liouville's formula \cite{MR2682403}. Thereby, the normal variational equations (NVE) along $\Gamma_0$ are
\begin{align}\label{eq3}
	&\xi''_j=\left(\varepsilon u_1^{2k-2}\left(t\right)-\mu_j\right)\xi_j,\quad j=2,\ldots,m.
\end{align}

Inspired
by  Yoshida \cite{MR923886}, we introduce the following finite branched covering map
\begin{equation}\label{eq17}
	\begin{split}
	&\overline{\Gamma}_0\rightarrow \mathbf{P}^1,\\
	&t\longmapsto z=\dfrac{\varepsilon}{k\mu_1}u_1^{2k-2}\left(t\right),
	\end{split}
\end{equation}
where $\overline{\Gamma}_0$ is the compact Riemann surface of the curve $v_1^2=\varepsilon u_1^{2k}/k-\mu_1u_1^2$ and $\mathbf{P}^1$ is the Riemann sphere.

Performing the Yoshida transformation \eqref{eq17}, the NVE \eqref{eq3} can be written as the hypergeometric differential equations in the new independent variable $z$
\begin{align}
	&\dfrac{d^2\xi_j}{dz^2}+\left(\frac{1}{z}+\frac{1}{2 (z-1)}\right)\dfrac{d\xi_j}{dz}-\left(\frac{\mu _j}{4 \mu _1(k-1)^2 z^2}+\frac{k \mu _1-\mu _j}{4\mu _1 (k-1)^2z (z-1)}\right)\xi_j=0,\tag{$\text{ANVE}_j$}\\
	&j=2,\ldots,m.\notag
\end{align}
The above differential system of equations is called the algebraic normal variational equations (ANVE, for short), and is denoted as
\begin{align}\label{eq4}
	&\text{ANVE}=\text{ANVE}_2\oplus\text{ANVE}_3\oplus\cdots\oplus\text{ANVE}_m.
\end{align}
Essentially, equation \eqref{eq4} is a direct sum in the more intrinsic sense of linear connections, see Chapter $2$ of \cite{MR1713573}  for more details.

From Theorem \ref{th4}, it follows that the identity components of the Galois groups of the NVE \eqref{eq3} and the ANVE \eqref{eq4} coincide. Obviously, the ANVE \eqref{eq4}  is integrable if and only if each $\text{ANVE}_j$ is integrable for $j=1,\ldots,m$. More precisely,  the identity component of the
Galois group of the ANVE is solvable if and only if the identity
 component of the Galois group of each $\text{ANVE}_j$ is solvable  for $j=1,\ldots,m$.

Now, we consider the $\text{ANVE}_2$:
\begin{align}\label{eq19}
&\dfrac{d^2\xi_2}{dz^2}+\left(\frac{1}{z}+\frac{1}{2 (z-1)}\right)\dfrac{d\xi_2}{dz}-\left(\frac{\mu _2}{4 \mu _1(k-1)^2 z^2}+\frac{k \mu _1-\mu _2}{4\mu _1 (k-1)^2z (z-1)}\right)\xi_2=0
\end{align}
with three singular points at $z=0,1,\infty$. Comparing  \eqref{eq19} with the general form of the
hypergeometric equation \eqref{hy} (see Appendix \ref{A}),  one can see that the exponents of \eqref{eq19} at singular points must
fulfill the following relations
\begin{align*}
&\alpha+\tilde{\alpha}=0,\quad \alpha\tilde{\alpha}=-\frac{\mu _2}{4 \mu _1(k-1)^2},\\
&\beta+\tilde{\beta}=\dfrac{1}{2},\quad \beta\tilde{\beta}=-\frac{k}{4(k-1)^2},\\
&\gamma+\tilde{\gamma}=\dfrac{1}{2},\quad \gamma\tilde{\gamma}=0.
\end{align*}
Thus, all the possibilities of the differences of exponents are
\begin{align}\label{eq21}
&\varrho=\pm\dfrac{1}{k-1}\sqrt{\dfrac{\mu_2}{\mu_1}},\tau=\pm\dfrac{1}{2}\left(1+\dfrac{2}{k-1}\right)\;\text{and}\;\varsigma=\pm\dfrac{1}{2}.
\end{align}
Let\begin{footnotesize}
$$\mathcal{S}=\left\{\pm\dfrac{1}{k-1}\left(1+\sqrt{\dfrac{\mu_2}{\mu_1}}\right),\pm\dfrac{1}{k-1}\left(1-\sqrt{\dfrac{\mu_2}{\mu_1}}\right),\pm\left(1+\dfrac{1}{k-1}\left(1+\sqrt{\dfrac{\mu_2}{\mu_1}}\right)\right),\pm\left(1+\dfrac{1}{k-1}\left(1-\sqrt{\dfrac{\mu_2}{\mu_1}}\right)\right)\right\}.$$
\end{footnotesize}Straightforward computations show that
\begin{align}\label{eq56}
&\varrho+\tau+\varsigma\in\mathcal{S},\;-\varrho+\tau+\varsigma\in\mathcal{S},\;\varrho-\tau+\varsigma\in\mathcal{S}\;\text{and}\;\varrho+\tau-\varsigma\in\mathcal{S}.
\end{align}

If equation \eqref{eq19} satisfies statement (i) of  Theorem \ref{A1}, by equation \eqref{eq56}, then
$$\text{either}\;\dfrac{1}{k-1}\left(\sqrt{\dfrac{\mu_2}{\mu_1}}+1\right)\;\text{or}\;\dfrac{1}{k-1}\left(\sqrt{\dfrac{\mu_2}{\mu_1}}-1\right)$$
must be an integer, that is,
\begin{align}\label{eq20}
&\dfrac{\mu_2}{\mu_1}\in\left\{\left(\left(k-1\right)\ell\pm1\right)^2\big|\ell\in \mathbb{N}\right\}.
\end{align}

If the statement (ii) of  Theorem \ref{A1} is fulfilled for equation \eqref{eq19}, from equation \eqref{eq21}, then its differences of exponents only conform the first row of Table \ref{Ta3}. Note that $k\geq3$. Therefore,
$$\pm\dfrac{1}{k-1}\sqrt{\dfrac{\mu_2}{\mu_1}}=\dfrac{1}{2}+\ell,\quad \ell\in\mathbb{Z},$$
that is,
$$\dfrac{\mu_2}{\mu_1}\in\left\{\dfrac{\left(k-1\right)^2\left(2\ell+1\right)^2}{4}\Bigg|\ell\in \mathbb{Z}\right\}.$$

Based on the discussion above,  the parameters $\mu_1$ and $\mu_2$ must satisfy
\begin{align}\label{eq22}
&\dfrac{\mu_2}{\mu_1}\in\left\{\left(\left(k-1\right)\ell\pm1\right)^2\big|\ell\in \mathbb{N}\right\}\bigcup\left\{\dfrac{\left(k-1\right)^2\left(2\ell+1\right)^2}{4}\Bigg|\ell\in \mathbb{Z}\right\}
\end{align}
if the identity components of the Galois groups of the NVE \eqref{eq3} is Abelian.

On the second invariant manifold $\mathcal{N}_2$, system \eqref{eq38} is written as
\begin{align}\label{eq6}
	&u'_2=v_2,\quad v'_2=-\mu_2u_2+\varepsilon u_2^{2k-1}
\end{align}
with Hamiltonian
\begin{align}\label{eq23}
&\tilde{h}=\dfrac{1}{2}v_2^2+\dfrac{1}{2}\mu_2u_2^2-\dfrac{\varepsilon}{2k}u_2^{2k}.
\end{align}
To solve equation \eqref{eq23}, we obtain
\begin{align}\label{eq48}
	&\dfrac{du_2}{dt}=\pm\sqrt{2\tilde{h}+\dfrac{\varepsilon}{k}u_2^{2k}-\mu_2u_2^2}.
\end{align}

Let $\widetilde{\Theta}\left(\tilde{h}\right)\in\mathbb{C}^2$  be a solution  of system \eqref{eq6} on the energy level $h$. So,
\begin{align}\label{eq58}
	&\widetilde{\varGamma}_{\tilde{h}}:=\left\{\left(0,0,u_2\left(t\right),v_2\left(t\right),0,\ldots,0\right)\in\mathbb{C}^{2m}\mid \left(u_2\left(t\right),v_2\left(t\right)\right)\in\widetilde{\Theta}\left(\tilde{h}\right)\right\}
\end{align}
is a particular solution of system \eqref{eq38}. Equation \eqref{eq48} has three equilibrium points $u_2=0$ and $u_2=\pm\sqrt[2k-2]{\mu_2k/\varepsilon}$ on the energy level $\tilde{h}=0$. In the same way as particular solution $\Gamma_0$, we can find a non-constant solution $u_2\left(t\right)$, which presents a non-equilibrium  particular solution $\widetilde{\Gamma}_0\in\widetilde{\varGamma}_0$.

Let $\bm{\eta}:=\left(\eta_1,\ldots,\eta_m\right)^T$ and $\bm{\tilde{\eta}}:=\left(\tilde{\eta}_1,\ldots,\tilde{\eta}_m\right)^T$. The variational equations (VE) along $\widetilde{\Gamma}_0$ is given by
\begin{align}\label{eq7}
	&\left(
	\begin{array}{c}
		\bm{\eta}'\vspace{1ex}\\
		\bm{\tilde{\eta}}'\\
	\end{array}\right)
	=
	\left(
	\begin{array}{cc}
		\mathbf{0}&\mathbf{I}\vspace{2ex}\\
		\mathbf{\tilde{\Lambda}}&\mathbf{0}\\
	\end{array}
	\right)
	\left(
	\begin{array}{c}
		\bm{\eta}\vspace{1ex}\\
	\bm{\tilde{\eta}}\\
	\end{array}\right),
\end{align}
where
$$\mathbf{\tilde{\Lambda}}:=\text{diag}\;\Big(\varepsilon u_2^{2k-2}\left(t\right)-\mu_1,\varepsilon\left(2k-1\right)u_2^{2k-2}\left(t\right)-\mu_2,\varepsilon u_2^{2k-2}\left(t\right)-\mu_3,\ldots, \varepsilon u_2^{2k-2}\left(t\right)-\mu_m\Big).$$
The VE \eqref{eq7} is also composed of  $m$ uncoupled Schr\"{o}dinger equations
\begin{align*}
	&\bm{\eta}''=\mathbf{\tilde{\Lambda}}\bm{\eta},
\end{align*}
that is,
\begin{align}
	&\eta''_1=\left(\varepsilon u_2^{2k-2}\left(t\right)-\mu_1\right)\eta_1,\notag\\
	&\eta''_2=\left(\varepsilon\left(2k-1\right)u_2^{2k-2}\left(t\right)-\mu_2\right)\eta_2,\label{eq64}\\
	&\eta''_j=\left(\varepsilon u_2^{2k-2}\left(t\right)-\mu_j\right)\eta_j,\quad j=3,\ldots,m.\notag
\end{align}
Using Liouville's formula \cite{MR2682403}, the second equation of \eqref{eq64} is solvable due to the fact that it has a solution $\eta_2=u'_2\left(t\right)$. Therefore, the corresponding normal variational equations ($\widetilde{\text{NVE}}$) along $\widetilde{\Gamma}_0$ are given by
\begin{equation}\label{eq8}
	\begin{split}
	&\eta''_1=\left(\varepsilon u_2^{2k-2}\left(t\right)-\mu_1\right)\eta_1,\\
&\eta''_j=\left(\varepsilon u_2^{2k-2}\left(t\right)-\mu_j\right)\eta_j,\quad j=3,\ldots,m.
	\end{split}
\end{equation}

Similarly, we can carry out the following Yoshida transformation
$$t\longmapsto z=\dfrac{\varepsilon}{k\mu_2}u_2^{2k-2}\left(t\right),$$
and transform $\widetilde{\text{NVE}}$ \eqref{eq8} into the algebraic normal variational equations ($\widetilde{\text{ANVE}}$):
\begin{align}
	&\dfrac{d^2\eta_1}{dz^2}+\left(\frac{1}{z}+\frac{1
	}{2 (z-1)}\right)\dfrac{d\eta_1}{dz}-\left(\frac{\mu _1}{4 \mu _2(k-1)^2 z^2}+\frac{k \mu _2-\mu _1}{4\mu _2 (k-1)^2z (z-1)}\right)\eta_1=0,\tag{$\widetilde{\text{ANVE}}_1$}\\
	&\dfrac{d^2\eta_j}{dz^2}+\left(\frac{1}{z}+\frac{1}{2 (z-1)}\right)\dfrac{d\eta_j}{dz}-\left(\frac{\mu _j}{4 \mu _2(k-1)^2 z^2}+\frac{k \mu _2-\mu _j}{4\mu _2 (k-1)^2z (z-1)}\right)\eta_j=0,\tag{$\widetilde{\text{ANVE}}_j$}\\
	&j=3,\ldots,m.\notag
\end{align}
The direct sum form of $\widetilde{\text{ANVE}}$ is
$\widetilde{\text{ANVE}}=\widetilde{\text{ANVE}}_1\oplus\widetilde{\text{ANVE}}_3\oplus\widetilde{\text{ANVE}}_4\oplus\cdots\oplus\widetilde{\text{ANVE}}_m$. For the $\widetilde{\text{ANVE}}_1$,  all the possibilities of the differences of exponents are
$$\varrho=\pm\dfrac{1}{k-1}\sqrt{\dfrac{\mu_1}{\mu_2}},\tau=\pm\dfrac{1}{2}\left(1+\dfrac{2}{k-1}\right)\;\text{and}\;\varsigma=\pm\dfrac{1}{2}.$$

 By the same discussions as NVE \eqref{eq3}, we obtain  that  the parameters $\mu_1$ and $\mu_2$ must satisfy
 \begin{align}\label{eq24}
& \dfrac{\mu_1}{\mu_2}\in\left\{\left(\left(k-1\right)\ell\pm1\right)^2\big|\ell\in \mathbb{N}\right\}\bigcup\left\{\dfrac{\left(k-1\right)^2\left(2\ell+1\right)^2}{4}\Bigg|\ell\in \mathbb{Z}\right\} 	
 \end{align}
if the identity components of the Galois groups of the $\widetilde{\text{NVE}}$ \eqref{eq8} is Abelian.

The conditions \eqref{eq22} and  \eqref{eq24} imply that
$$\dfrac{\mu_2}{\mu_1}\geq1\;\text{and}\;\dfrac{\mu_1}{\mu_2}\geq1,$$
respectively. This
contradicts our assumption $\mu_1\neq\mu_2$. Consequently, either the identity components of the Galois groups of the NVE \eqref{eq3} or  $\widetilde{\text{NVE}}$ \eqref{eq8} is not Abelian.
By Theorem \ref{th5}, the lemma follows.

The proof is finished.
\qed

\begin{lem}\label{le3}
If $k\geq3$, $\mu_1\neq\mu_2$ and $\mu_1\mu_2=0$, then Hamiltonian system \eqref{eq38} is not meromorphic integrable.
\end{lem}

\Proof Our proof will be distinguished two cases: either $ \mu_1=0, \mu_2\neq0$, or $\mu_1\neq0, \mu_2=0$.

Case 1: $\mu_1=0$ and $\mu_2\neq0$. For this case, we also restrict system \eqref{eq38} on the invariant manifold $\mathcal{N}_1$. Namely,
\begin{align}\label{eq15}
	&u'_1=v_1,\quad v'_1=\varepsilon u_1^{2k-1}
\end{align}
with Hamiltonian
\begin{align}\label{eq25}
	&h=\dfrac{1}{2}v_1^2-\dfrac{\varepsilon}{2k}u_1^{2k}.
\end{align}

Analogously, we also consider the particular solution $\Gamma_0$ in the proof of Lemma \ref{le2}, and compute the $\widehat{\text{NVE}}$ along $\Gamma_0$
\begin{equation}\label{eq45}
	\begin{split}
&\xi''_j=\left(\varepsilon u_1^{2k-2}\left(t\right)-\mu_j\right)\xi_j,\quad j=2,\ldots,m.
	\end{split}
\end{equation}

Doing the change of variable
$$t\longmapsto z=\dfrac{\varepsilon}{2\mu_2}u_1^{2k-2}\left(t\right),$$
we attain the algebraic normal variational equations ($\widehat{\text{ANVE}}$):
\begin{align}
	&\dfrac{d^2\xi_2}{d z^2}+\dfrac{3}{2z}\dfrac{d\xi_2}{d z}-\frac{k\left(2z-1\right)}{8 (k-1)^2z^3}\xi_2=0,\tag{$\widehat{\text{ANVE}}_2$}\\
	&\dfrac{d^2\xi_j}{d z^2}+\dfrac{3}{2z}\dfrac{d\xi_j}{d z}-\frac{k \left(2\mu _2z-\mu _j\right)}{8(k-1)^2\mu _2 z^3}\xi_j=0,\quad j=3,\ldots, m, \tag{$\widehat{\text{ANVE}}_j$}
\end{align}
and denote by
\begin{align}\label{eq5}
	&\widehat{\text{ANVE}}=\widehat{\text{ANVE}}_2\oplus\widehat{\text{ANVE}}_3\oplus\cdots\oplus\widehat{\text{ANVE}}_m.
\end{align}
Performing the classical transformation (see \eqref{eq63})
$$\xi_2=\chi\exp\left(-\dfrac{3}{4}\int \dfrac{dz}{z} \right)=\chi z^{-3/4},$$
the $\widehat{\text{ANVE}}_2$ becomes
\begin{align}\label{eq26}
&\chi''=r\left(z\right)\chi,
\end{align}
where
\begin{align}\label{eq27}
&r\left(z\right)=-\left(\dfrac{(k-3) (3 k-1)}{16 (k-1)^2 z^2}+\dfrac{ k}{8(k-1)^2 z^3}\right).
\end{align}
Then, the poles of  $r\left(z\right)$ is $z=0$. The order of  $z=0$ and $z=\infty$ is $o\left(0\right)=3$ and $o\left(\infty\right)=2$, respectively.  Using Proposition \ref{pr1} to equation  \eqref{eq26}, only types (ii) or (iv) of Theorem \ref{th6} can appear. Working the second part of Kovacic's algorithm (see Appendix \ref{B}), we obtain that
$$\mathcal{E}_\infty=\left\{2+\ell\sqrt{1-\dfrac{(k-3) (3 k-1)}{4(k-1)^2 }}\;\Bigg|\;\ell=0,\pm2\right\}\bigcap\mathbb{Z}=
\begin{cases}
	\left\{0,2,4\right\},\;\text{if}\;k=3,\\
	\left\{2\right\},\;\text{if}\;k\geq4.
\end{cases}\;\text{and}\;\mathcal{E}_0=\left\{3\right\}.$$
Straightforward computations show that the  number $d=d\left(\bm{\varpi}\right)=\left(\varpi_\infty-\varpi_0\right)/2$
is not an integer. Therefore, type  (iv) of Theorem \ref{th6}  holds. This means that the identity component of the
Galois group of the $\widehat{\text{ANVE}}$ \eqref{eq5} is not Abelian. Thereby, the identity component of the
Galois group of the  $\widehat{\text{NVE}}$ \eqref{eq45} is also not Abelian. From Theorem \ref{th5}, it follows that the Hamiltonian system \eqref{eq38} for $k\geq3$ is meromorphic non-integrable with $\mu_1=0$ and $\mu_2\neq0$.

 Case 2: $\mu_1\neq0$ and $\mu_2=0$. Substituting $\mu_2=0$ into \eqref{eq8}, we get the normal variational equations along $\widetilde{\Gamma}_0$:
 \begin{equation}\label{eq46}
 	\begin{split}
&\eta''_1=\left(\varepsilon u_2^{2k-2}\left(t\right)-\mu_1\right)\eta_1,\\
&\eta''_j=\left(\varepsilon u_2^{2k-2}\left(t\right)-\mu_j\right)\eta_j,\quad j=3,\ldots,m.
 	\end{split}
 \end{equation}

After the change of variable
$$t\longmapsto z=\dfrac{\varepsilon}{2\mu_1}u_2^{2k-2}\left(t\right),$$
equations \eqref{eq46} become the algebraic normal variational equations
 \begin{align*}
 	&\dfrac{d^2\eta_1}{d z^2}+\dfrac{3}{2z}\dfrac{d\eta_1}{d z}-\frac{k\left(2z-1\right)}{8 (k-1)^2z^3}\eta_1=0,\\
 	&\dfrac{d^2\eta_j}{d z^2}+\dfrac{3}{2z}\dfrac{d\eta_j}{d z}-\frac{k \left(2\mu _1z-\mu _j\right)}{8(k-1)^2\mu _1 z^3}\eta_j=0,\quad j=3,\ldots, m.
 \end{align*}
  The analysis is exactly the same as  Case 1.  Thus, the Hamiltonian system \eqref{eq38} for $k\geq3$ is meromorphic non-integrable with $\mu_1\neq0$ and $\mu_2=0$.

 This lemma holds.
\qed

Now we can state the main results of this paper.
\begin{thm}\label{th1}
	The Hamiltonian system \eqref{eq38} is meromorphic completely integrable if and only if $k=2$ and $\left(\mu_1,\ldots,\mu_m\right)\in\mathbb{R}^m$, or $k\geq3$ and $\mu_1=\mu_2=\cdots=\mu_m$.
\end{thm}

\Proof The sufficiency was given by Theorem 3.4 of \cite{TZ}.

For the necessity, it is enough to prove that the Hamiltonian system \eqref{eq38} is meromorphic non-integrable if $k$ and $\left(\mu_1,\ldots,\mu_m\right)$ are not satisfied the sufficient condition of Theorem \ref{th1}. This means that $k\geq3$ and there exists a positive integer $j_0\in\left\{2,\ldots,m\right\}$ such that $\mu_1\neq\mu_{j_0}$. We can assume without loss of generality that $j_0=2$, that is, $\mu_1\neq\mu_2$, because in the other case  one can interchange respectively the roles of $\mu_{j_0}$ and $\mu_2$, and $u_{j_0}$ and $u_2$. By Lemma \ref{le2} and Lemma \ref{le3}, we get the necessity.

The proof is ended.
\qed

\section*{Acknowledgments}

The first author is partially supported by the National Natural Science Foundation of China (Grant no. 12401205 and Grant no. 12371182), the Guangdong Basic and Applied Basic Research Foundation (Grant no. 2023A1515110430), the Science and Technology Planning Project of Guangzhou (Grant no. 2025A04J3404), the China Postdoctoral Science Foundation (Grant no. 2025M773097), the Postdoctoral Fellowship Program (Grade C) of China Postdoctoral Science Foundation (Grant no. GZC20230970) and the Fundamental Research Funds for the Central Universities (Grant no. 21624347). The second author is partially supported by the National Natural Science Foundation of China (Grant no. 11790273). 

%

\section*{Author Declarations}
\subsection*{Data availability statements}
The authors declare that no datasets were generated or analysed in this study.
\subsection*{Conflict of Interest}
The authors declare that they have no conflict of interest.

\section*{Appendix}
\setcounter{equation}{0}
\setcounter{subsection}{0}
\setcounter{figure}{0}
\setcounter{table}{0}
\renewcommand{\theequation}{\thesection.\arabic{equation}}
\renewcommand{\thesubsection}{\thesection.\arabic{subsection}}
\renewcommand{\thetable}{\thesection.\arabic{table}}
\renewcommand{\thefigure}{\thesection.\arabic{figure}}
\begin{appendices}

	\section{Hypergeometric equation}\label{A}
	The hypergeometric equation is a  second order differential equation over the Riemann sphere $\mathbf{P}^1$ with three  regular singular points \cite{MR2682403,MR0178117}. Let us consider the following form of hypergeometric equation with three singular points at $z=0,1,\infty$
	\begin{align}\label{hy}
		&\dfrac{d^2\zeta}{dz^2}+\left(\dfrac{1-\alpha-\tilde{\alpha}}{z}+\dfrac{1-\gamma-\tilde{\gamma}}{z-1}\right)\dfrac{d\zeta}{dz}+\left(\dfrac{\alpha\tilde{\alpha}}{z^2}+\dfrac{\gamma\tilde{\gamma}}{\left(z-1\right)^2}+\dfrac{\beta\tilde{\beta}-\alpha\tilde{\alpha}-\gamma\tilde{\gamma}}{z\left(z-1\right)}\right)\zeta=0,
	\end{align}
	where $\left(\alpha,\tilde{\alpha}\right)$, $\left(\gamma,\tilde{\gamma}\right)$ and $\left(\beta,\tilde{\beta}\right)$ are the exponents at the respective singular points, and meet the Fuchs relation
	$\alpha+\tilde{\alpha}+\gamma+\tilde{\gamma}+\beta+\tilde{\beta}=1$.
	The exponent differences can be defined as
	$\varrho=\alpha-\tilde{\alpha}$, $\varsigma=\gamma-\tilde{\gamma}$ and $\tau=\beta-\tilde{\beta}$.
	
	The following theorem goes back to Kimura \cite{MR277789}, whose gave
	necessary and sufficient
	conditions for solvability of the identity component of the differential Galois group of  \eqref{hy}.
	
	\begin{table}[H]
		\centering
		\caption{Schwarz table with $l,s,\upsilon\in\mathbb{Z}.$}\vspace{2ex}
		\begin{tabular}{|c|c|c|c|c|}
			\hline
			$1$&$1/2+l$&$1/2+s$&Arbitrary complex number&\xrowht{6pt}\\
			\hline
			$2$&$1/2+l$&$1/3+s$&$1/3+\upsilon$&\xrowht{6pt}\\
			\hline
			$3$&$2/3+l$&$1/3+s$&$1/3+\upsilon$&$l+s+\upsilon$ even\xrowht{6pt}\\
			\hline
			$4$&$1/2+l$&$1/3+s$&$1/4+\upsilon$&\xrowht{6pt}\\
			\hline
			$5$&$2/3+l$&$1/4+s$&$1/4+\upsilon$&$l+s+\upsilon$ even\xrowht{6pt}\\
			\hline
			$6$&$1/2+l$&$1/3+s$&$1/5+\upsilon$&\xrowht{6pt}\\
			\hline
			$7$&$2/5+l$&$1/3+s$&$1/3+\upsilon$&$l+s+\upsilon$ even\xrowht{6pt}\\
			\hline
			$8$&$2/3+l$&$1/5+s$&$1/5+\upsilon$&$l+s+\upsilon$ even\xrowht{6pt}\\
			\hline
			$9$&$1/2+l$&$2/5+s$&$1/5+\upsilon$&$l+s+\upsilon$ even\xrowht{6pt}\\
			\hline
			$10$&$3/5+l$&$1/3+s$&$1/5+\upsilon$&$l+s+\upsilon$ even\xrowht{6pt}\\
			\hline
			$11$&$2/5+l$&$2/5+s$&$2/5+\upsilon$&$l+s+\upsilon$ even\xrowht{6pt}\\
			\hline
			$12$&$2/3+l$&$1/3+s$&$1/5+\upsilon$&$l+s+\upsilon$ even\xrowht{6pt}\\
			\hline	
			$13$&$4/5+l$&$1/5+s$&$1/5+\upsilon$&$l+s+\upsilon$ even\xrowht{6pt}\\
			\hline
			$14$&$1/2+l$&$2/5+s$&$1/3+\upsilon$&$l+s+\upsilon$ even\xrowht{6pt}\\
			\hline		
			$15$&$3/5+l$&$2/5+s$&$1/3+\upsilon$&$l+s+\upsilon$ even\xrowht{6pt}\\
			\hline							
		\end{tabular}%
		\label{Ta3}%
	\end{table}%
	\begin{thm}\label{A1}{\rm (\cite{MR277789})} The identity component of the Galois group of the hypergeometric equation \eqref{hy} is solvable if and only if  either
		\begin{itemize}
			\item [\emph{(i)}] at least one of the four numbers    $\varrho+\tau+\varsigma,-\varrho+\tau+\varsigma,\varrho-\tau+\varsigma,\varrho+\tau-\varsigma$ is an odd integer, or
			\item [\emph{(ii)}] the numbers $\varrho$ or $-\varrho$, $\varsigma$ or $-\varsigma$ and $\tau$ or  $-\tau$ belong (in an arbitrary order)
			to some of the following fifteen families, see Table \ref{Ta3}.
		\end{itemize}
	\end{thm}

\section{Kovacic's algorithm}\label{B}
Let $\mathbb{C}\left(z\right)$ be the field of rational functions  in the variable $z$ with complex coefficients. Consider the second order linear differential equation
\begin{align}\label{eq28}
	&\chi''=r\left(z\right)\chi,\quad r\left(z\right) \in \mathbb{C}\left(z\right).
\end{align}
It is well known that  the differential Galois group $G$ of equation \eqref{eq28} is an algebraic subgroup of $\text{SL}\left(2,\mathbb{C}\right)$.  In 1986, Kovacic \cite{MR839134}
characterized all possible types of $G$ as follows.

\begin{thm}\label{th6}{\rm (\cite{MR839134})}
	The differential Galois group $G$ of equation \eqref{eq28} is conjugated to one of the following four types:
	\begin{itemize}
		\item [\emph{(i)}] $G$ is conjugated to a subgroup of a triangular group, and equation \eqref{eq28}  admits a solution of the form $\chi=\exp\left(\int \omega \right)$ with $\omega \in\mathbb{C}\left(z\right)$.
		\item [\emph{(ii)}] $G$  is conjugate to a subgroup of
		\begin{align*}
			&\mathcal{G}=\left\{
			\left(
			\begin{array}{cc}
				\mathfrak{a}&0\\
				0&\mathfrak{a}^{-1}\\
			\end{array}
			\right)
			\Bigg|\mathfrak{a}\in \mathbb{C}\setminus\left\{0\right\}\right\}
			\bigcup
			\left\{
			\left(
			\begin{array}{cc}
				0&\mathfrak{a}\\
				\mathfrak{a}^{-1}&0\\
			\end{array}
			\right)
			\Bigg|\mathfrak{a}\in \mathbb{C}\setminus\left\{0\right\}
			\right\},
		\end{align*}
		and equation \eqref{eq28}  admits a solution of the form $\chi=\exp\left(\int \omega \right)$, where $\omega$ is algebraic of degree $2$ over $\mathbb{C}\left(z\right)$.
		\item [\emph{(iii)}] $G$ is finite and all solutions of equation \eqref{eq28} are algebraic over $\mathbb{C}\left(z\right)$.
		\item [\emph{(iv)}] $G=\emph{SL}\left(2,\mathbb{C}\right)$ and equation \eqref{eq28} does not admit a Liouvillian solution.
	\end{itemize}
\end{thm}

Let $r\left(z\right)=\mathfrak{p}\left(z\right)/\mathfrak{q}\left(z\right)$
with $\mathfrak{p}\left(z\right),\mathfrak{q}\left(z\right)\in\mathbb{C}\left[z\right]$ relatively prime. The \emph{pole} of $r\left(z\right)$ is a zero of $\mathfrak{q}\left(z\right)$ and \emph{the order of the pole} is the multiplicity of the zero of $\mathfrak{q}\left(z\right)$.  \emph{The order of $r\left(z\right)$ at $\infty$} is defined by $\deg \mathfrak{q}-\deg \mathfrak{p}$.
Kovacic \cite{MR839134} also
provided the necessary conditions for types (i), (ii), or (iii) in  Theorem \ref{th6} to occur.
\begin{prop}\label{pr1}{\rm (\cite{MR839134})}
	For the first three types in Theorem \ref{th6}, the necessary conditions of occurrence are respectively as follows:
	\begin{description}
		\item[\emph{\bf Type (i)}] Each pole of $r\left(z\right)$ must have even order or else have order $1$. The order of $r\left(z\right)$ at $\infty$ must be
		even or else be greater than $2$.
		\item[\emph{\bf Type (ii)}] The rational function $r\left(z\right)$ must have at least one pole that either has odd order greater than $2$ or else has
		order $2$.
		\item[\emph{\bf Type (iii)}] The order of a pole of $r\left(z\right)$ cannot exceed $2$ and the order of $r\left(z\right)$ at $\infty$ must be at
		least $2$. If the partial fraction decomposition of $r\left(z\right)$ is
		$$r\left(z\right)=\sum_i\dfrac{\alpha_i}{\left(z-c_i\right)^2}+\sum_j\dfrac{\beta_j}{z-b_j},$$
		then $\sqrt{1+4\alpha_i}\in\mathbb{Q}$ for each $i$, $\sum_j\beta_j=0$, and if $\Delta=\sum_i\alpha_i+\sum_j\beta_j$, then $\sqrt{1+4\Delta}\in\mathbb{Q}$.
	\end{description}
\end{prop}

\begin{rem}\label{re1}
	The general second order linear differential equation
	\begin{align*}
		&y''=a_1y'+a_2,\quad a_1,a_2\in \mathbb{C}\left(z\right),
	\end{align*}
	can be transformed into the form \eqref{eq28} with
	$$r\left(z\right)=\dfrac{a_1^2}{4}-\dfrac{a'_1}{2}+a_2$$
	via the change of the variable
	\begin{align}\label{eq63}
		&y=\exp\left(\dfrac{1}{2}\int a_1 dz\right)\chi.
	\end{align}
\end{rem}

Here, we recall  the second part of  Kovacic's algorithm \cite{MR839134}. Let $r\left(z\right)\in C\left(z\right)$ and $\Upsilon$ be the set of poles of $r\left(z\right)$. Set $\chi''=r\chi$.
	\begin{itemize}
		\item [\bf{Step 1.}] To each pole $c\in\Upsilon$, we calculate the set  $\mathcal{E}_c$ as follows.
		\begin{itemize}
			\item [(i)] If the pole $c$ is of order $1$, then $\mathcal{E}_c=\left\{4\right\}$.
			\item [(ii)] If the pole $c$ is of order $2$ and $b$ is  the coefficient of $1/\left(z-c\right)^2$ in the partial fraction decomposition of $r\left(z\right)$, then
			\begin{align}\label{eq29}
				&\mathcal{E}_c=\left\{2+\ell\sqrt{1+4b}\;\big|\;\ell=0,\pm2\right\}\bigcap\mathbb{Z}.
			\end{align}
			\item [(iii)] If the pole $c$ is of order $o\left(c\right)>2$, then $\mathcal{E}_c=\left\{o\left(c\right)\right\}$.
			\item [(iv)] If the order of $r$ at $\infty$ is $o\left(\infty\right)>2$, then $\mathcal{E}_c=\left\{0,2,4\right\}$.
			\item [(v)] If the order of $r$ at $\infty$ is $2$ and $b$ is  the coefficient of $1/z^2$ in the Laurent expansion  of $r\left(z\right)$ at $\infty$, then
			\begin{align}\label{eq42}
				&\mathcal{E}_c=\left\{2+\ell\sqrt{1+4b}\;\big|\;\ell=0,\pm2\right\}\bigcap\mathbb{Z}.
			\end{align}
			\item [(vi)] If the order of $r$ at $\infty$ is $o\left(\infty\right)<2$, then $\mathcal{E}_c=\left\{o\left(\infty\right)\right\}$.
		\end{itemize}
		\item [\bf{Step 2.}] Let  $\bm{\varpi}=\left(\varpi_c\right)_{c\in\Upsilon}$ be a element in the Cartesian product $\prod_{c\in\Upsilon}\mathcal{E}_c$ with $\varpi_c\in\mathcal{E}_c$. Define number
		\begin{align}\label{eq43}
			&d:=d\left(\bm{\varpi}\right)=\dfrac{1}{2}\left(\varpi_\infty-\sum_{c\in\Upsilon}\varpi_c\right).
		\end{align}
		We try to find all elements  $\bm{\varpi}$ such that $d$ is a non-negative integer, and retain such elements to perform Step 3.  If there is no such element  $\bm{\varpi}$,  then statement (ii) of Theorem \ref{th6} is impossible.
		\item [\bf{Step 3.}] For each $\bm{\varpi}$ retained from Step 2,  we introduce the rational function
		$$\theta=\dfrac{1}{2}\sum_{c\in\Upsilon}\dfrac{\varpi_c}{z-c}.$$
		Then, we seek a monic polynomial $P$ of degree $d$ defined in \eqref{eq43} such that
		\begin{align*}
			&P'''+3\theta P''+\left(3\theta^2+3\theta'-4r\right)P'+\left(\theta''+3\theta\theta'+\theta^3-4r\theta-2r'\right)P=0,
		\end{align*}
where monic polynomial $P$ is a polynomial with the leading coefficent $1$.
If such polynomial $P$ does not exist for all  elements $\bm{\varpi}$ retained from Step 2,  then statement (ii) of Theorem \ref{th6} is untenable.

Assume that  such a polynomial $P$ exists. Let $\phi=\theta+P'/P$ and $\omega$ be a root of
$$\omega^2-\phi \omega+\left(\dfrac{1}{2}\phi'+\dfrac{1}{2}\phi^2-r\right)=0.$$
Then, $\chi=\exp\left(\int \omega \right)$ is a solution of differential equation $\chi''=r\chi$.
	\end{itemize}
  \end{appendices}

\end{document}